\documentclass[10pt]{amsart}
\title[
]{
The behaviour of curvature 
functions \\ at cusps and inflection points}

\date{February 22, 2010.}
\usepackage[usenames]{color}

\usepackage[dvips]{graphicx} 
\usepackage[mathscr]{eucal}
\usepackage{mathrsfs}
\usepackage{amscd}
\usepackage{amsfonts}
\usepackage{amsmath,amsthm,amssymb}
\usepackage{latexsym}
\usepackage{bm}
\usepackage[dvips]{graphics}
\numberwithin{equation}{section}

\usepackage{amsthm}
\theoremstyle{plain}
 \newtheorem{theorem}{Theorem}[section]
 \newtheorem*{theorem*}{Theorem}

 \newtheorem*{lemma*}{Lemma}
 \newtheorem{proposition}[theorem]{Proposition}
 
 \newtheorem*{fact*}{Fact}

 \newtheorem{lemma}[theorem]{Lemma}
 \newtheorem{corollary}[theorem]{Corollary}
\theoremstyle{remark}
 \newtheorem{definition}[theorem]{Definition}
 \newtheorem{remark}[theorem]{Remark}
 \newtheorem*{remark*}{Remark}
 
 \newtheorem{example}[theorem]{Example}
\numberwithin{equation}{section}

\newcommand{\R}{\boldsymbol{R}}

\renewcommand{\phi}{\varphi}
\renewcommand{\epsilon}{\varepsilon}
\newcommand{\op}{\operatorname}

\newcommand{\mb}[1]{{\mathbf #1}}
\newcommand{\pmt}[1]{{\begin{pmatrix} #1  \end{pmatrix}}}
\newcommand{\dy}{\displaystyle}

\author{Shohei~Shiba}
\address[Shiba]{%
   Department of Mathematics, Graduate School of Science,
   Osaka University,
   Toyonaka, Osaka 560-0043,
   Japan
}
\email{}

\author{Masaaki Umehara}
\address[Umehara]{%
   Department of Mathematics, Graduate School of Science,
   Osaka University,
   Toyonaka, Osaka 560-0043,
   Japan
}
\email{umehara@math.sci.osaka-u.ac.jp}


\thanks{
The second author was partially 
supported by the Grant-in-Aid for 
Scientific Research (A) No.22244006, 
Japan Society for the Promotion of Science.}

\begin{document}

\begin{abstract}
At a $3/2$-cusp of a given plane curve $\gamma(t)$,
both of the Euclidean curvature $\kappa_g$ and
the affine curvature $\kappa_A$  diverge.
In this paper, we show that each of $\sqrt{|s_g|}\kappa_g$
and $(s_A)^2 \kappa_A$ (called the {\it Euclidean 
{\em and} affine normalized curvature}, respectively) 
at a $3/2$-cusp is a $C^\infty$-function of the variable $t$,
where $s_g$ (resp. $s_A$)
is  the Euclidean (resp. affine) arclength parameter 
of the curve corresponding to the $3/2$-cusp
$s_g=0$ (resp. $s_A=0$).
Moreover, we give a 
characterization of the behaviour of 
the curvature functions $\kappa_g$ and $\kappa_A$
at $3/2$-cusps.
On the other hand, inflection points 
are also singular points of curves in affine geometry.
We give a similar characterization of 
affine curvature functions near generic inflection points.
As an application, new affine invariants
of $3/2$-cusps and generic inflection points are given. 
\end{abstract}
\maketitle

\section{Introduction.}
Let $\gamma(t)$ be a smooth (i.e. $C^\infty$) curve
in the plane $\R^2$
defined on an open interval containing $t=0$.
The origin $t=0$ is called a {\it singular point} 
of $\gamma$ if $\dot \gamma(t)=d\gamma(t)/dt$ 
vanishes at $t=0$.
Moreover, a singular point $t=0$ is called
a {\it $3/2$-cusp} if there exist a suitable coordinate 
change $t=t(s)$ and a local diffeomorphism $\Phi$
of $\R^2$ at $\gamma(0)$ such that
$\Phi\circ \gamma\circ s(t)={}^t\!(t^2, t^3)$, where
\lq$t$\rq\ denotes the transpose operation on matrices.
It is well-known that a singular point $t=0$
is a $3/2$-cusp if and only if
$[\ddot\gamma(0),\gamma^{(3)}(0)]\ne 0$ holds,
where $[\mb a,\mb b]$ denotes the determinant, that is
$$
[\mb a,\mb b]=a_1b_2-a_2b_1, \qquad \mb a=\pmt{a_1\\ a_2},\,\,
\mb b=\pmt{b_1\\ b_2}.
$$
Moreover, a $3/2$-cusp $t=0$ is 
called a {\it positive cusp}
(resp. a {\it negative cusp})
if $[\ddot\gamma(0),\gamma^{(3)}(0)]$ is positive
(resp. negative).
This signature of cusps is invariant 
under an orientation preserving
diffeomorphism of $\R^2$. 
If one reverses the orientation of the curve,
the signature of the cusp changes.
An invariant called the (Euclidean)
{\it cuspidal curvature}
$\mu^{}_g$ at $3/2$-cusps is given as follows
(which was introduced in \cite{U},
and its fundamental properties are given in \cite{SUY})
\begin{equation}\label{eq:m_g}
\mu^{}_g:=\frac{[\ddot\gamma(0),
\gamma^{(3)}(0)]}{|\ddot \gamma(0)|^{5/2}},
\end{equation}
which is independent of 
orientation preserving isometries of $\R^2$. 
The sign of the cuspidal curvature coincides with that of
cusps. Let $\gamma_1(t)$ and $\gamma_2(t)$ be two $3/2$-cusps
at $t=0$. Suppose that $\gamma_1$ and $\gamma_2$
have the same {\it cuspidal curvature}, then there exists  
orientation preserving isometry $T$
of $\R^2$ and a parametrization $u=u(t)$ near $t=0$
such that $|T\circ \gamma_2\circ u(t)-\gamma_1(t)|$
has order higher than $t^3$.

If we denote by $s_g$ the (Euclidean) arclength
parameter of the curve such that $s_g=0$ corresponds to
a $3/2$-cusp, then we show firstly in this paper that 
$\sqrt{|s_g|}\, \kappa_g$ 
(called the \lq  normalized 
curvature function\rq)
is a $C^\infty$-function of the variable $t$,
which induces the cuspidal curvature as follows:
\begin{theorem} \label{thm:main_g}
Let $\gamma(t):(-\epsilon,\epsilon)$ 
$(\epsilon>0)$ be a $3/2$-cusp in 
the Euclidean plane $\R^2$
such that $t=0$ is a $3/2$-cusp.
Then $\tau=\op{sgn}(t)\sqrt{|s_g|}$ can be
taken to be a coordinate of $\gamma$ at $t=0$
{$($called the \em  half-arclength parameter)}
and $\sqrt{|s_g|}\,\kappa_g$ is a $C^\infty$-function 
of $t$ $($and $\tau)$. Moreover,
it holds that
\begin{equation}\label{eq:lim_g0}
\lim_{t\to 0}\sqrt{|s_g|}\, \kappa_g=\frac{\mu_g}{2\sqrt{2}},
\end{equation}
where $s_g=\int_0^t |\dot \gamma(u)|du$ is the
$($Euclidean$)$ arclength parameter and $\mu_g$ is the
cuspidal curvature of $\gamma$ at $t=0$ given in 
\eqref{eq:m_g}.
Conversely, if we take a $C^\infty$-function
$f(\tau)$ such that $f(0)\ne 0$,
then there exists a $3/2$-cusp
such that
\begin{equation}\label{eq:N_g}
\sqrt{|s_g|}\,\kappa_g=f(\tau),
\end{equation}
and $\tau$ is the half-arclength parameter.
\end{theorem}
The formula \eqref{eq:lim_g0} was not given
in \cite{U0} and \cite{SUY}.
For each smooth function $\phi(s)$, 
it is well-known that there exists a regular 
curve $\gamma(s)$ with arclength parameter 
whose curvature function  is $\phi(s)$. 
The last assertion of the theorem 
is an analogue of this fact for $3/2$-cusps. 
Later, we also show that the same
assertion holds for $3/2$-cusps in an arbitrarily
given Riemannian 2-manifold (cf. Theorem \ref{thm:g_g}).

On the other hand, in affine geometry, 
Izumiya-Sano \cite{IS} pointed out 
the fact that the affine evolute 
having $3/2$-cusps corresponding to 
sextactic points is exactly analogous 
to the fact that the  Euclidean evolute 
has cusps corresponding to vertices.
Related to this work, Giblin and Sapiro \cite{GS}
studied the affine distance symmetry set
from the viewpoint of singularity theory. 
These two works suggest that the $3/2$-cusp is 
also an important object in affine geometry.
In this paper, we define a new affine invariant 
for $3/2$-cusps called the
{\it affine cuspidal curvature}, by
\begin{equation}\label{eq:mu_A}
\mu_A:=\left.
\frac{24[\ddot \gamma,\gamma^{(3)}]\,[\ddot \gamma,\gamma^{(5)}]
+60[\ddot \gamma,\gamma^{(3)}]\,[\gamma^{(3)},\gamma^{(4)}]
-35[\ddot \gamma,\gamma^{(4)}]^2}
{[\ddot \gamma,\gamma^{(3)}]^{12/5}}\right |_{t=0},
\end{equation}
where $t=0$ is a $3/2$-cusp of $\gamma(t)$.
It is invariant under equi-affine transformations
and independent of the choice of an 
orientation of the curve $\gamma(t)$. 
(An {\it equi-affine transformation} is an affine transformation
whose Jacobian is $\pm 1$ with respect to the canonical 
coordinate system of $\R^2$.)
Let $m,n$ be two mutually prime integers.
Here (and also throughout in this paper), we use the convention 
for the fractional order of exponent
as follows
\begin{equation}\label{eq:exponent}
t^{m/n}:=(-1)^{mn}|t|^{m/n}.
\end{equation}
For example, $t^{1/2}$ is equal to $\sqrt{|t|}$.
As an analogue of Euclidean cuspidal curvature,
we show the following two assertions:

\begin{theorem}\label{thm:main}
Let $\gamma_1(t)$ and $\gamma_2(t)$ be two 
 $3/2$-cusps in the affine plane $\R^2$
at $t=0$. Then $\gamma_1$ and $\gamma_2$
have the same affine cuspidal curvature
if and only if
there exists an equi-affine
transformation $T$ 
of $\R^2$ and a parametrization $u=u(t)$ near $t=0$
such that $du/dt>0$ and
$|T\circ \gamma_2\circ u(t)-\gamma_1(t)|$
has order higher than $t^5$.
\end{theorem}

\begin{theorem}\label{thm:k_A}
Let $\gamma(t):(-\epsilon,\epsilon)$ 
$(\epsilon>0)$ be 
a $3/2$-cusp
in the affine plane $\R^2$
such that $t=0$ is a $3/2$-cusp
and $s_A$ is the affine arclength parameter
$($cf. \eqref{eq:S_A}$)$.
Then 
$\tau=(s_A)^{3/5}$ {\rm (cf. \eqref{eq:exponent})
can be taken to be a coordinate of $\gamma$ at $t=0$
{$($called the \em  $3/5$-arclength parameter$)$},
and $f:=(s_A)^2\kappa_A$ is a $C^\infty$-function 
of $t$ $($and $\tau)$,
called the {\em  normalized affine curvature 
function} 
(see \eqref{eq:def_k} for definition of $\kappa_A$).
Moreover,
it satisfies $f(0)=4/25$, $\dot f(0)=0$ 
and
\begin{equation}\label{eq:lim_A2}
\lim_{t\to 0} \frac{(s_A)^2\kappa_A-{4}/{25}}{\tau^2}
=\frac{1}{220}\sqrt[5]{\frac{20}3}\mu_A.
\end{equation}}
Conversely, if we take a $C^\infty$-function
$f(\tau)$ such that $f(0)=4/25$
and $\dot f(0)=0$, 
then there exists a $3/2$-cusp
whose normalized affine curvature function is
$f(\tau)$ with respect to the $3/5$-length
parameter.
\end{theorem}

Later, we generalize Theorem \ref{thm:k_A}
for $3/2$-cusps in an arbitrarily
given 2-manifold with equi-affine structure
(cf. Theorem \ref{thm:g_A}).

A point $t=0$ on a regular curve $\gamma(t)$
is called an {\it inflection point} if
it satisfies 
$[\dot \gamma(0), \ddot \gamma(0)]=0$.
An inflection point is called {\it generic}
if $[\dot \gamma(0), \gamma^{(3)}(0)]\ne 0$ holds. 
Affine geometry of curves usually requires 
higher order derivatives than those in 
Euclidean geometry, and this fact is often 
connected to several interesting phenomena 
different from Euclidean geometry: 
For example, inflection points are singular points 
in affine geometry, as well as $3/2$-cusps.
In Section 3, we show an analogue of Theorems 
\ref{thm:main} and \ref{thm:k_A} for inflection points.
In particular, the normalized (affine) 
curvature function $f(t)=(s_A)^2\kappa_A$
is a $C^\infty$-function at a generic inflection point $t=0$,
which satisfies $f(0)=-5/16$ and also 
the following non-trivial identity
(see Theorem \ref{thm:k_I} and Theorem \ref{thm:g_I})
\begin{equation}\label{eq:hidden2}
-\frac{9}{7}\left(\frac{\left([\dot \gamma(0),\gamma^{(4)}(0)]
+[\ddot \gamma(0),\gamma^{(3)}(0)]\right)}{[\dot \gamma(0),
\gamma^{(3)}(0)]}\right)\dot f(0)
+32 \dot f(0)^2+9 \ddot f(0)=0.
\end{equation}
This corresponds to the fact that 
$\dot f(0)$ vanishes at $3/2$-cusps, 
as in Theorem \ref{thm:k_A}.

\section{The natural equation of cusps in Euclidean Geometry}

We denote by $C^\infty_0(\R)$ 
the set of germs of real-valued $C^\infty$-functions defined at $t=0$.
Let $\gamma(t)$ be a curve defined on an interval $(-\delta,\delta)$
for $\delta>0$. We suppose that $t=0$ is a $3/2$-cusp, 
namely it satisfies $\dot \gamma(0)=\mathbf 0$ and
$[\ddot \gamma(0),\gamma^{(3)}(0)]\ne 0$.
Then it can be easily checked that the Euclidean
curvature function 
\begin{equation}\label{eq:k_g}
\kappa_g(t):=\frac{[\dot \gamma(t),\ddot 
\gamma(t)]}{|\dot \gamma(t)|^3}
\end{equation}
diverges. More precisely, the following assertion
holds:

\begin{lemma}\label{prop:k_g}
Let $t=0$ be a $3/2$-cusp of the curve $\gamma(t)$
in the Euclidean plane $\R^2$.
Then, $\tau:=\op{sgn}(t)\sqrt{|s_g|}$ 
can be taken as a local coordinate of the
curve $\gamma$ at $t=0$.
Moreover, $\sqrt{|s_g|}\kappa_g$ belongs to $C^\infty_0(\R)$, 
and \eqref{eq:lim_g0} holds.
\end{lemma}

\begin{proof}
Without loss of generality, 
we may assume that $\gamma(0)=\mathbf 0$. 
We may set
\begin{equation}\label{eq:g00}
\gamma(t)=\frac{\ddot\gamma(0)}{2}t^2
+\frac{\gamma^{(3)}(0)}{6}t^3 +t^4 \Gamma(t),
\end{equation}
where $\Gamma(t)$ is a $\R^2$-valued 
$C^\infty$-function.
Then we get the following expression
\begin{equation}\label{eq:g1}
[\dot\gamma(t),\ddot\gamma(t)]=
t^2\left(\frac{[\ddot\gamma(0),\gamma^{(3)}(0)]}{2} 
+t \phi_1(t)\right),
\end{equation}
where $\phi_1(t)\in C^\infty_0(\R)$. Similarly,
there exists
$\phi_2(t)\in C^\infty_0(\R)$ such that
\begin{equation}\label{eq:g2}
|\dot\gamma(t)|=|t|\,\left| \ddot\gamma(0)+t\phi_2(t)\right|.
\end{equation}
Then \eqref{eq:g1} and \eqref{eq:g2} imply that
$|t|\kappa_g(t)\in C^\infty_0(\R)$.
On the other hand, applying Lemma \ref{lem:fractional} 
below by setting $\alpha=1$
and $\phi(t):=\left| \ddot\gamma(0)+t\phi_2(t)\right|$,
we can conclude that 
$$
f(t):=\frac{s_g(t)}{\op{sgn}(t)t^2}=
\frac{|s_g(t)|}{t^2}=
\left(\frac{\sqrt{|s_g(t)|}}{|t|}\right)^2
$$
is a $C^\infty$-function such that $f(0)>0$.
Then 
$$
\tau:=\op{sgn}(t)\sqrt{|s_g(t)|}=t\sqrt{f(t)}
$$
is a $C^\infty$-function of $t$.
Since $d\tau(0)/dt=f(0)>0$,
the function $\tau$ can be taken as a local coordinate of the
curve $\gamma$ at $t=0$.
Since $|t|\kappa_g\in C^\infty_0(\R)$,
the function
\begin{equation}\label{eq:product}
\sqrt{|s_g^{}|}\kappa_g=(|t|\kappa_g)\sqrt{f(t)}
\end{equation}
is also a $C^\infty$-function of $t$. Finally,
the formula \eqref{eq:lim_g0} follows directly from
\eqref{eq:g1}, \eqref{eq:g2}, \eqref{eq:product}
and $f(0)=|\ddot \gamma(0)|/2$ (cf. (A.1) 
in the appendix).
\end{proof}

\begin{remark} \label{rmk:1/2}
A parametrization $t$
of the $3/2$-cusp $\gamma(t)$ in $\R^2$
is the half-arclength parameter
if and only if it satisfies $|\dot \gamma|=2|t|$.
\end{remark}

\medskip
\noindent
{\it Proof of Theorem \ref{thm:main_g}.}\\
The first parts of the assertion have been proved 
in Lemma \ref{prop:k_g}.
So it is sufficient to show the last assertion.
We take a $C^\infty$-function $f(\tau)$.
Let $\gamma(\tau)$ be a curve defined by
\begin{equation}
\label{eq:integral}
\gamma(\tau):=2 \int_0^\tau u \pmt{\cos \theta(u) \\ \sin \theta(u)}du
\qquad \left(\theta(\tau):=2\int_0^\tau f(u)du \right).
\end{equation}
Then it holds that $|\dot \gamma|=2\tau$, which implies that
$\tau$ is the half-arclength parameter.
Moreover, one can directly check that $\gamma(\tau)$
has a $3/2$-cusp satisfying \eqref{eq:N_g}.
\qed

\begin{example}
The arclength parameter of
the cusp 
$$
\gamma(t)=a\pmt{t^2 \\ t^3}
\qquad  (a>0)
$$
is given by
$$
s_g=\op{sgn}(t)(\phi(t)-8)a
\qquad
\left(\phi(t):=\left(9 t^2+4\right)^{3/2}\right),
$$
which is a continuous function, but not smooth at $t=0$.
The half-arclength parameter $\tau=\op{sgn}(t)a\sqrt{\phi(t)-8}$
is a smooth function at $t=0$ and
the normalized curvature function
is 
$$
\frac{2 (\phi(t)-8)^{1/2}}{\sqrt{3a} t \phi(t)}
=\frac{3}{4\sqrt{a}}-\frac{297 t^2}{128\sqrt{a}}+o\left(t^3\right)
\qquad
\left(\phi(t):=\left(9 t^2+4\right)^{3/2}\right),
$$
which implies that the cuspidal curvature at $t=0$
is equal to $3/\sqrt{2a}$.
\end{example}

\begin{example}
The half-arclength parameter of 
the cycloid
\begin{equation}\label{eq:cyc}
a\pmt{t-\sin t\\-1+\cos t}\qquad (a>0)
\end{equation}
at $t=0$ is equal to $2\sqrt{2a}\sin(t/4)$
and the normalized curvature function is given by
$$
\frac{1}{2 \sqrt{2a}\cos(t/4)}=
\frac{1}{2 \sqrt{2a}}+
\frac{t^2}{64 \sqrt{2a}}+o\left(t^3\right),
$$
which implies that the cuspidal curvature
at $t=0$ is equal to $1/\sqrt{a}$.
\end{example}

\begin{figure}[htb]
 \begin{center}
        \includegraphics[height=4.0cm]{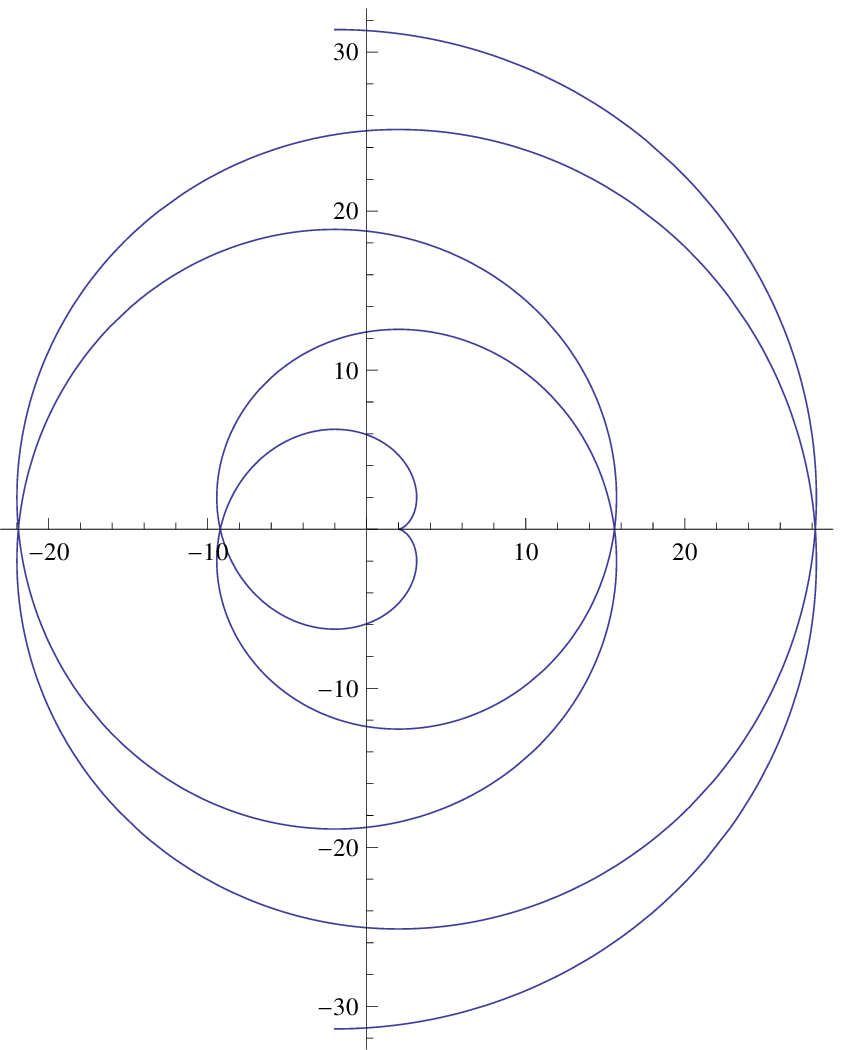}
\caption{The canonical $3/2$-cusp ($a=1$) in Euclidean geometry
}\label{fig:canonical}
\end{center}
\end{figure}

\begin{example}
The curve given by (cf. Figure \ref{fig:canonical})
$$
\frac1{2a^2}\pmt{{2 a \tau \sin (2 a \tau)+
\cos (2 a \tau)}\\ \sin (2 a \tau)-
2 a \tau \cos (2 a \tau)}
$$
is called the {\it canonical $3/2$-cusp}, which
has the property that $\sqrt{|s_g|}\kappa_g$ is 
identically equal to $a$, where $\tau$ is the
half-arclength parameter.
\end{example}

To finish this section, we generalize 
Theorem \ref{thm:main_g} to $3/2$-cusps
in an arbitrary Riemannian 2-manifold:
Let $\gamma(t)$ ($|t|<\delta$) be a regular 
curve in a given oriented Riemannian $2$-manifold $(M^2,g)$.
We denote by $\Omega_g$ the unit area element
of $M^2$. By definition,
$\Omega_g(\mb e_1,\mb e_2)=1$ holds for a positively
oriented orthonormal frame $\mb e_1,\mb e_2$ of $(M^2,g)$.
For the sake of simplicity, we set
$$
[v,w]:=\Omega_g(v,w)\qquad (v,w\in T_pM^2,\,\, p\in M^2).
$$
This notation fits with the previous one,
since $\Omega_g(v,w)=\op{det}(v,w)$ holds on 
the Euclidean plane $\R^2$. 
Then the geodesic curvature $\kappa_g$ of the curve
is defined by exactly the same
formula \eqref{eq:k_g}, where
$$
\ddot \gamma(t)=\nabla_t\dot\gamma(t),\qquad
|\dot \gamma(t)|=\sqrt{g(\dot \gamma(t),\dot \gamma(t))}.
$$
On the other hand, let
$\gamma(t)$ ($|t|<\delta$) be a $3/2$-cusp 
at $t=0$ in $M^2$. Then we define the {\it cuspidal
curvature} $\mu_g$ at $t=0$ by \eqref{eq:m_g}
like as in the case of $\R^2$, where
$$
\gamma^{(3)}=\nabla_t\ddot\gamma,\quad
\gamma^{(4)}=\nabla_t\gamma^{(3)}.
$$
Then the following assertion holds, which is a 
generalization of Theorem \ref{thm:main_g}.

\begin{theorem}\label{thm:g_g}
Let $\gamma(t):(-\delta,\delta)$ 
$(\delta>0)$ be a curve in
an oriented Riemannian $2$-manifold $(M^2,g)$
such that $t=0$ is a $3/2$-cusp.
Then the same assertion as in Theorem 
\ref{thm:main_g} holds.
\end{theorem}

To prove the assertion, we prepare the following assertion:

\begin{lemma}\label{lem:parallel}
Let $X(t)$ be a $C^\infty$-vector field along $\gamma$ such that
\begin{equation}\label{eq:diffn}
X(0)=\dot X(0)=\cdots =X^{(n)}(0)=\bf 0,
\end{equation}
where $\dot X=\nabla_t X$ and $X^{(i)}:=\nabla_tX^{(i-1)}$
$(i=1,2,3,...)$.
Then there exist $\delta>0$ and a 
$C^\infty$-vector field $Y(t)$ along 
$\gamma(t)$ $(0\le t<\delta)$
such that $X(t)=t^{n+1}Y(t)$.
\end{lemma}

\begin{proof}
Let $(\mb e_1(t),\mb e_2(t))$ be a parallel frame field along $\gamma$.
We may set
$$
X(t):=f_1(t)\mb e_1(t)+f_2(t)\mb e_2(t),
$$
and then it holds that
$$
X^{(i)}(t)=f_1^{(i)}(t)\mb e_1(t)+f_2^{(i)}(t)\mb e_2(t)
\qquad (i=1,2,...,n).
$$
Then \eqref{eq:diffn} is equivalent to the condition
$f_j(0)=f^{(1)}_j(0)\cdots=f_j^{(n)}(0)=0$ ($j=1,2$).
Then, there exists a $C^\infty$-function
$g_j(t)\in C_0^\infty(\R)$ such that
$f_j(t)=t^{n+1}g_j(t)$ ($j=1,2$) and
$X=t^{n+1}(g_1 \mb e_1+g_2\mb e_2)$ holds.
\end{proof}

\noindent
\medskip
({\it Proof of Theorem \ref{thm:g_g}}.)
The assertion of Lemma \ref{prop:k_g} holds for
$3/2$-cusps, however, the proof should be modified:
Let $v$ be a tangent vector of $M^2$ at $\gamma(0)$.
Then we denote by $P_v(t)$ the parallel vector field 
along $\gamma(t)$ such that $P_v(0)=v$.
We set
\begin{equation}
X(t):=\dot\gamma(t)-P_2
t-\frac12 P_{3}t^2,
\end{equation}
where $P_2:=P_{\ddot \gamma(0)}(t)$ 
and $P_3:=P_{\gamma^{(3)}(0)}(t)$.
Then it holds that
$$
X(0)=\dot X(0)=\ddot X(0)=\mb 0.
$$
By Lemma \ref{lem:parallel},
there exists a vector field $\Gamma(t)$
along $\gamma(t)$ such that
$$
\dot\gamma(t)=P_2 t+\frac12 P_3 t^2+
t^3 \Gamma(t).
$$
Then it holds that
$$
\ddot\gamma(t)=P_2+P_{3} t+
t^2 (3\Gamma(t)+t \dot \Gamma(t)).
$$
Using these expressions, one can easily prove 
the same assertion of
Lemma \ref{prop:k_g} for a curve $\gamma$ 
on the Riemannian manifold $(M^2,g)$.
In particular,
$f(t):=\sqrt{|s_g|}\kappa_g$
is a smooth function of $t$, and 
$\tau:=\op{sgn}(t)\sqrt{|s_g|}$ gives a
local coordinate of the curve at the $3/2$-cusp.
We can also check the identity 
\eqref{eq:lim_g0} by the completely same argument
as in the proof of Lemma \ref{prop:k_g}.

However, our previous proof of
the last assertion
of Theorem \ref{thm:main_g} 
cannot be applied, since \eqref{eq:integral}
holds only for the case that the ambient space is $\R^2$.
To prove it here, we need the following new idea: 
Let $\gamma(\tau)$ be a $3/2$-cusp at $\tau=0$, 
and suppose that $\tau$ is the half-arclength parameter.
We use the notations 
\begin{equation}\label{eq:dot}
\dot \gamma:=\frac{d \gamma}{d\tau},\quad
\ddot \gamma:=\frac{d^2 \gamma}{d\tau^2},\quad
\gamma^{(3)} :=\frac{d^3 \gamma}{d\tau^3},\quad
\gamma^{(4)}:=\frac{d^4 \gamma}{d\tau^4},\,\,\cdots,
\end{equation}
and 
\begin{equation}\label{eq:prime}
\gamma':=\frac{d \gamma}{ds},\quad
\gamma'':=\frac{d^2 \gamma}{ds^2},\quad
\gamma^{[3]}:=\frac{d^3 \gamma}{ds^3},\quad
\gamma^{[4]}:=\frac{d^4 \gamma}{ds^4},\,\,\cdots,
\end{equation}
where $s=s_g$. Since $s_g=\tau^2$,
the following identities hold:
\begin{equation}\label{eq:g12}
\dot \gamma=2\tau\gamma', \qquad
\ddot \gamma=
2 \mb e + 4 \tau f(\tau)\mb n,
\end{equation}
where $\mb e=\gamma'$ and $\mb n$ is the 
unit normal vector field along 
$\gamma$ such that $[\mb e,\mb n]=1$. 
Here, we used the following identities
\begin{equation}\label{eq:frenet2}
\dot{\mb e}=2\tau \mb e= 2\tau \kappa_g\mb n=
2f(\tau)\mb n,\qquad
\dot{\mb n}=2\tau \mb n= -2f(\tau)\mb n.
\end{equation}
Using \eqref{eq:g12}, one can reprove
the identity \eqref{eq:lim_g0}, since both sides of
\eqref{eq:lim_g0} are independent of the choice of the
parameter $t$.
On the other hand, we have that
$$
\gamma^{(3)}
=-8 \tau f^2 \mb e+4 (\tau \dot f+2f) \mb n.
$$
Since $\ddot \gamma$ and $\gamma^{(3)}$ at $\tau=0$
are linearly independent, using 
the Schmidt orthonormalization of 
the frame $\{\ddot \gamma(\tau), \gamma^{(3)}(\tau)\}$,
we get the following $C^\infty$-orthonormal frame
field along $\gamma$
$$
\mb u_1:=\frac1{\sqrt{1+4\tau^2 f^2}}(\mb e+2\tau f \mb n),\qquad
\mb u_2:=\frac1{\sqrt{1+4\tau^2 f^2}}(2\tau f \mb e-\mb n).
$$
Then it holds that
\begin{equation}\label{eq:ode_g}
\dot \gamma=\frac{2\tau}{\sqrt{1+4\tau^2 f^2}}
(\mb u_1-2\tau f \mb u_2).
\end{equation}
Applying \eqref{eq:frenet2}, we have that
\begin{equation}\label{eq:ode_g2}
\frac{d}{d\tau}(\mb u_1,\mb u_2)=
(\mb u_1,\mb u_2)\frac{2(2f+4\tau^2 f^3+\tau \dot f)}{1+4\tau^2 f^2}\pmt{0 &  
-1 \\ 1 & 0}.
\end{equation}
Using this observation, we shall now prove the last assertion.
We fix a $C^\infty$-function $f\in C^\infty_0(\R)$
satisfying $f(0)\ne 0$.
Let $\gamma(\tau)$ be a solution 
of the system of
ordinary differential equations \eqref{eq:ode_g}
and \eqref{eq:ode_g2} with the initial conditions
$$
\gamma(0)=\mb 0,\quad \mb u_1(0)=\pmt{1 \\ 0},\quad \mb u_2(0)=\pmt{0 \\ 1}.
$$
By \eqref{eq:ode_g}, we have
\begin{equation}\label{eq:abs}
|\dot \gamma|=2|\tau|.
\end{equation}
On the other hand, one can verify that
\begin{align}\label{eq:gamma1}
&\ddot \gamma=2\sqrt{1+4\tau^2 f^2} \mb u_1, \\
\label{eq:gamma2}
&\gamma^{(3)}(0)=8 f(0) \mb u_2(0),
\end{align}
and then we have that
$$
[\ddot \gamma(0),\gamma^{(3)}(0)]=16 f(0)\ne 0,
$$
which implies that $\gamma(\tau)$ has a $3/2$-cusp
at $\tau=0$. 
Then \eqref{eq:abs} and Remark \ref{rmk:1/2} yield that
$\tau$ is the half-arclength parameter.
By \eqref{eq:ode_g} and \eqref{eq:gamma1},
the identity \eqref{eq:N_g} follows immediately.
\qed

We denote by $S^*_0(\R,M^2)$ the 
set of germs of $C^\infty$-curves $\gamma(t)$ 
on $M^2$ which give a $3/2$-cusp at $t=0$.
Then the map 
\begin{equation}\label{eq:map_g}
\mathcal F_g:S^*_0(\R,M^2)
\ni \gamma(t)\mapsto \sqrt{|s_g(t)|}\kappa_g(t)
\in C^\infty_0(\R)
\end{equation}
is defined, namely $\mathcal F_{g}(\gamma)$ is the
normalized curvature function of the $3/2$-cusp $\gamma$.

\begin{corollary}\label{cor:image_g}
The image of the map $\mathcal F_g$ coincides with 
the subset
$$
\Sigma_g:=\left\{f\in C^\infty_0(\R)\,;\, f(0)\ne 0
\right\}.
$$
\end{corollary}

\begin{proof}
Obviously the image of $\mathcal F_g$ is contained
in $\Sigma_g$.
By applying Theorem \ref{thm:g_g},
for each $f\in \Sigma_g$, 
there exist $\delta>0$ and
a $3/2$-cusp
$\gamma:(-\delta,\delta)\to (M^2,g)$
such that $\sqrt{|s_g(t)|}\kappa_g(t)=f(t)$
and $t$ is the half-arclength parameter.
Namely, $\mathcal F_g(\gamma)=f$ holds. 
\end{proof}

\section{Affine geometry of cusps}

Let $\gamma(t)$ be a curve in $\R^2$
defined on an interval $(-\delta,\delta)$
for $\delta>0$. We suppose that $t=0$ is a $3/2$-cusp.
The formula \eqref{eq:lim_g0} suggests how to construct new 
invariants on singular points, and we shall now introduce
an affine invariant of $3/2$-cusps in the same manner.
In this section, we assume that $\gamma(t)$ has no inflection
points for $t\in (-\delta,\delta)$.
It is classically known that the affine curvature
function is defined by (cf. \cite{FB} and \cite{U})
\begin{equation}\label{eq:def_k}
\kappa_A=
\frac{
3[\dot \gamma,\ddot \gamma]\,
[\dot \gamma,\gamma^{(4)}]
+
12[\dot \gamma,\ddot \gamma]\,[\ddot \gamma,\gamma^{(3)}]
-5[\dot \gamma,\gamma^{(3)}]^2
}{9[\dot \gamma,\ddot \gamma]^{8/3}},
\end{equation}
which is invariant under equi-affine transformations and 
is independent of the choice of a parametrization $t$
and of the orientation of $\gamma$.
Here $[\dot \gamma,\ddot \gamma]^{8/3}$ is positive
because of our convention \eqref{eq:exponent}.
The {\it affine arclength} of $\gamma(t)$ is defined by
\begin{equation}\label{eq:S_A}
s_A^{}(t)
:=\int_{0}^t 
\left|[\dot \gamma(u),\ddot \gamma(u)]^{1/3}\right| du,
\end{equation}
which is invariant under equi-affine transformations
in $\R^2$. This parameter satisfies 
$$
\epsilon_A:=[\dot \gamma(s_A^{}),\ddot \gamma(s_A^{})]=\pm 1,
$$
and 
$$
\kappa_A^{}(s_A)=
\epsilon_A [\ddot \gamma(s_A^{}),\gamma^{(3)}(s_A^{})].
$$
It can be easily checked that 
$\kappa_A$ diverges at $3/2$-cusps. 
To prove Theorem \ref{thm:main},
it is sufficient to 
show the following assertion, since 
orientation reversing equi-affine transformation
of $\R^2$ preserve $\mu_A$ but reverse the signature of
$3/2$-cusps:

\begin{proposition}\label{prop:equi1}
Let $\gamma(t)$ be a positive $3/2$-cusp at $t=0$
in the affine plane $\R^2$.
Then there exist an orientation preserving equi-affine 
transformation $T$ and a coordinate 
change $t=t(u)$ such that $dt/du>0$ and
\begin{equation}
T\circ \gamma\circ t(u)=
\pmt{u^2\\ u^3 +\dy\frac{\mu_A u^5}{80\sqrt[5]{54}}}
+o(u^5),
\end{equation}
where $o(u^5)$ is a term of order higher than $u^5$.
\end{proposition}

\begin{proof}
We may assume that 
$[\ddot \gamma(0),\gamma^{(3)}(0)]=1$
by a suitable homothetic change of $t$.
By a parallel translation of $\R^2$, we may set
$\gamma(0)=\mb 0$.
If we set 
$$
T:=(\ddot\gamma(0),\gamma^{(3)}(0))\in \op{SL}(2,\R),
$$
then the new curve
$
\Gamma_0(t):=T^{-1}\gamma(t)
$
satisfies
$$
\Gamma_0(0)=\dot \Gamma_0(0)=\mathbf 0,\quad
\ddot\Gamma_0(0)=\pmt{1 \\ 0},\quad 
\Gamma^{(3)}_0(0)=\pmt{0 \\ 1}.
$$
So we can write 
$$
\Gamma_0(t)=\pmt{{t^2a_0(t)^2}/2 \\ {t^3b_0(t)}/6},
$$
where $a_0(t)$ and $b_0(t)$ are smooth functions at $t=0$
satisfying $a_0(0)=b_0(0)=1$.
If we set
$
s:=t a_0(t),
$
then it gives a new parametrization 
of $\Gamma_0$ such that 
$$
\Gamma_0(v)=\pmt{v^2/2\\ v^3b(v)/6},
$$
where $b(v)$ is a smooth function at $v=0$
satisfying $b(0)=1$.
Next, we set $v=12^{1/5}w$ and 
$$
\Gamma(w):=\pmt{k^{-1} & 0\\ 0 & k}\Gamma_0(w)\qquad (k:=\frac{6}{12^{3/5}}).
$$
Then we can write 
$$
\Gamma(u)=\pmt{w^2 \\ w^3(1+wB(w))},
$$
where $B(w)$ is a  smooth function at $w=0$.
Next, we set
$$
\hat\Gamma(w):=\pmt{1 & \xi \\
     0 & 1} \Gamma(w)
=
\pmt{w^2(1 + \xi w(1+wB(w)) \\ w^3(1+w B(w))},
$$
where $\xi$ is a constant.
If we take a new parameter 
$
u:=w\sqrt{1+\xi w+w^2 \xi B(w)},
$
then it satisfies
$$
u^3=w^3+\frac{3\xi  w^4}{2}+o(w^4).
$$
So if we set
$
\xi:=2B(0)/3,
$
it holds that
$$
\hat \Gamma(u)=\pmt{u^2\\ u^3+c u^5+o(u^5)}.
$$
By \eqref{eq:mu_A}, we 
have $c={\mu_A}/({80 \sqrt[5]{54}})$. 
\end{proof}

To prove Theorem \ref{thm:k_A},
we prepare the following assertion:

\begin{proposition}\label{prop:k_A}
Let $t=0$ be a $3/2$-cusp of a curve $\gamma(t)$
in the affine plane $\R^2$.
Then $\tau:=(s_A)^{3/5}$ can be taken as a new local 
parametrization of $\gamma(t)$ at $t=0$.
$($We shall call $\tau$ the 
{\em $3/5$-arclength parameter}.$)$
\end{proposition}

\begin{proof}
Without loss of generality, we may assume that
$\gamma(0)=\mathbf 0$. Then we have the expression
$\gamma(t)=t^2 \Gamma(t)$,
where $\Gamma(t)$ is a smooth $\R^2$-valued function
around $t=0$.
Using this, we can write
$
[\dot \gamma,\ddot \gamma]=t^2 a(t),
$
where $a\in C^\infty_0(\R)$ satisfies $a(0)\ne 0$.
Now applying Lemma \ref{lem:fractional}
by setting $\alpha=2/3$ and $\phi(t)=a(t)$,
the function
$
f(t)=t^{-5/3}s_A
$
is a $C^\infty$-function satisfying $f(0)\ne 0$.
Thus
$$
\tau:=(s_A)^{3/5}=t f(t)^{3/5}
$$
gives a new local parametrization of the curve $\gamma$.
\end{proof}

\begin{remark}\label{rmk:3/5}
Let $\gamma(\tau)$ be a $3/2$-cusp at $\tau=0$.
If $\tau$ is the $3/5$-arclength parameter,
it holds that
\begin{equation}\label{eq:tau-s}
\dot \gamma=\frac53\tau^{2/3}\gamma', \qquad
\ddot \gamma=
\frac{5 \left(2 \gamma'+5 \tau^{5/3}\gamma''\right)}{9 \tau^{1/3}},
\end{equation}
where we use the notations as in 
\eqref{eq:dot} and \eqref{eq:prime},
namely, the dot (resp. the prime)
means the derivative with respect to 
the $3/5$-arclength parameter
$\tau$ (resp. the affine arclength parameter $s_A$).
Since $[\gamma',\gamma'']=\pm 1$,
it can be easily seen that the parameter $t$ is the
$3/5$-arclength parameter if and only if
$[\dot \gamma(t),\ddot \gamma(t)]=\pm 125 t^2/27$.
\end{remark}

\begin{corollary}\label{cor:k_A}
Let $t=0$ be a $3/2$-cusp of a curve $\gamma(t)$
in $\R^2$.
Then $(s_A)^2 \kappa_A$ is a $C^\infty$-function at $t=0$.
$($As mentioned in the introduction, we call
$(s_A)^2 \kappa_A$ the {\em normalized affine curvature function}
at the cusp $t=0$.$)$
\end{corollary}

\begin{proof}
Without loss of generality, we may assume that
$t=0$ is a positive cusp.
Then it holds that $[\dot\gamma,\ddot \gamma]>0$.
We may assume that
$\gamma(0)=\mb 0$. Then we may set
\begin{equation}\label{eq:a0}
\gamma(t)=
\frac{\ddot\gamma(0)}{2!}t^2
+\frac{\gamma^{(3)}(0)}{3!}t^3
+t^4 \Gamma(t),
\end{equation}
where $\Gamma(t)$ is a smooth $\R^2$-valued function
around $t=0$.
For the sake of simplicity, we set
$
d_{23}:=[\ddot\gamma(0),\gamma^{(3)}(0)](\ne 0).
$
Then \eqref{eq:a0} yields the following expression
$$
[\dot \gamma,\ddot \gamma]=t^2 a_1(t),\quad
[\dot \gamma,\gamma^{(3)}]=t a_2(t),\quad
[\ddot \gamma,\gamma^{(3)}]=a_3(t),\quad
[\dot \gamma,\gamma^{(4)}]=t a_4(t),
$$
where $a_j(t)$ ($j=1,2,3,4$) are smooth functions of $t$.
In particular, $\kappa_A$ as in
\eqref{eq:def_k}
satisfies
$\kappa_A=t^{-10/3}\psi_1(t)$, 
where $\psi_1(t)\in C^\infty_0(\R)$.
Since $s_A$ can be expressed by
$s_A=t^{5/3}\psi_2(t)$ ($\psi_2\in C^\infty_0(\R)$)
(cf. Proposition \ref{prop:k_A}),
we get the assertion. 
\end{proof}

\medskip
\noindent
{\it Proof of Theorem \ref{thm:k_A}.}
Let $\gamma(\tau)$ be a $3/2$-cusp at $\tau=0$
and $\tau$ the $3/5$-arclength parameter.
Differentiating \eqref{eq:tau-s},
we have the following identity
\begin{equation}\label{eq:g3}
\gamma^{(3)}=\frac{5(-2 \gamma'+30 \gamma''
 \tau^{5/3}+25 \tau^{10/3}\gamma^{[3]} )}{27 \tau^{4/3}},
\end{equation}
where we use the notations as in 
\eqref{eq:dot} and \eqref{eq:prime}.
By Corollary \ref{cor:k_A}, we may set
$
(s_A)^{2}\kappa_A=f(\tau).
$
Then substituting the identity
$\gamma^{[3]}=-\tau^{-10/3}f\gamma'$ into \eqref{eq:g3},
we have
\begin{equation}\label{eq:g30}
\gamma^{(3)}
=
\frac5{{27 \tau^{4/3}}}\left(-(2+25f) \gamma' 
+30\tau^{5/3}\gamma'' \right).
\end{equation}
Differentiating \eqref{eq:g30} by 
using $\gamma^{[3]}=-\tau^{-10/3}f\gamma'$ again,
we have
\begin{equation}\label{eq:c4a}
\gamma^{(4)}=
-\frac{5}{81 \tau^{7/3}}
\left((75 \tau \dot f+50 f-8) \gamma '
+5 \tau^{5/3} (25 f-4) \gamma ''\right).
\end{equation}
Since
$[\gamma',\gamma'']=1$,
it holds that
\begin{equation}\label{eq:d14}
[\dot\gamma(\tau),\gamma^{(4)}(\tau)]
=
\frac{5^3}{3^5} \biggl(4 - 25 f(\tau)\biggr).
\end{equation}
Since $\dot\gamma(0)=\mb 0$, we have 
$f(0)={4}/{25}$.
So there exists a function $g\in C_0^\infty(\R)$
such that
$
f=({4}/{25})+\tau g.
$
Then it holds that
\begin{equation}\label{eq:d34}
[\gamma^{(3)},\gamma^{(4)}]
=
\frac{5^5 \left(18 \tau \dot g+
25 \tau g^2+36 g\right)}{3^7 \tau}.
\end{equation}
Since the left hand side is smooth at $\tau=0$,
we can conclude that $g(0)=0$.
So we may set
$
f=({4}/{25})+\tau^2 h,
$
where $h\in C^\infty_{0}(\R)$.
Using this expression
and the relation $[\gamma',\gamma'']=1$,
we have that
\begin{align}\label{eq:23a}
&[\dot\gamma,\ddot\gamma]
=\frac{125 \tau^2}{27}, \quad
[\dot\gamma,\gamma^{(3)}]=\frac{250 \tau}{27},\quad
[\ddot\gamma,\gamma^{(3)}]
=\frac{5^3}{3^5} \left(18+25 \tau^2 h\right),\\
&\label{eq:14a}
[\dot\gamma,\gamma^{(4)}]=-\frac{5^5}{3^5} \tau^2 h(t),\quad
[\ddot\gamma,\gamma^{(4)}]=\frac{5^5}{3^5} 
\tau \left(\tau \dot h+2 h\right),\\ 
\label{eq:34a}
&[\gamma^{(3)},\gamma^{(4)}]=
\frac{5^5}{3^7} \left(18 \tau 
\dot h+25 \tau^2 h^2+54 h\right).
\end{align}
Since
$$
\frac{5^5}{3^5} 
\tau \left(3\dot h+\tau \ddot h\right)\\ 
=
\frac{d}{d\tau}
[\ddot\gamma,\gamma^{(4)}]=
[\gamma^{(3)},\gamma^{(4)}]+[\ddot\gamma,\gamma^{(5)}],
$$
we have 
$
[\ddot\gamma(0),\gamma^{(5)}(0)]=-4h(0)(5/3)^5.
$
Using these relations
and \eqref{eq:mu_A}, 
one can easily check that
$$
\frac1{220}\sqrt[5]{\frac{20}{3}}\mu_A=h(0)
\left(=\lim_{\tau\to 0} \frac{(s_A)^2 \kappa_A-4/25}{\tau^2}\right).
$$
Since both sides of \eqref{eq:lim_A2}
are independent of the choice of parameters of
the curve, this proves the identity \eqref{eq:lim_A2}. 

Next we prove the last part of the theorem.
By \eqref{eq:23a}, $\gamma(\tau)$ must satisfy
the initial conditions 
\begin{equation}\label{eq:initial}
\dot \gamma(0)=\mathbf 0,\quad
[\ddot \gamma(0),\gamma^{(3)}(0)]=
\frac{250}{27}.
\end{equation}
Since $\ddot \gamma(0),\gamma^{(3)}(0)$ 
are linearly independent,
we can write
\begin{equation}\label{ode:A1}
\dot \gamma
=\alpha_1 \ddot\gamma+\alpha_2 \gamma^{(3)},
\quad
\gamma^{(4)}
=\beta_1 \ddot\gamma+\beta_2 \gamma^{(3)},
\end{equation}
where
\begin{equation}\label{eq:alpha}
\begin{split}
\alpha_1&=\frac{[\dot\gamma,\gamma^{(3)}]}{[\ddot\gamma,\gamma^{(3)}]}
=\frac{18\tau}{18+25 \tau^2 h},\qquad 
\alpha_2=-\frac{[\dot\gamma,\ddot\gamma]}{[\ddot\gamma,\gamma^{(3)}]}
=-\frac{9\tau^2}{18+25 \tau^2 h}, \\
\beta_1&
=-\frac{[\gamma^{(3)},\gamma^{(4)}]}{[\ddot\gamma,\gamma^{(3)}]}=
-\frac{25 \left(18 \tau \dot h+25 \tau^2 h^2+54 h\right)}
{9 \left(18+25 \tau^2 h\right)},\\
\beta_2&
=-\frac{[\ddot\gamma,\gamma^{(4)}]}{[\ddot\gamma,\gamma^{(3)}]}
=\frac{25\tau \left(\tau \dot h+2 h\right)}
{18+25\tau^2h}. 
\end{split}
\end{equation}
In particular, 
$\alpha_1,\alpha_2,\beta_1,\beta_2\in C^\infty_0(\R)$.
We now fix $h\in C^\infty_0(\R)$,
and take a solution $\gamma(\tau)$ 
of the ordinary differential 
equation
\begin{equation}\label{ode:A}
\dot\gamma=\alpha_1 \xi+\alpha_2 \eta,\qquad 
\dot\xi=\eta, \qquad \dot \eta=\beta_1 \xi+\beta_2 \eta
\end{equation}
with the initial conditions
\begin{equation}\label{eq:initial_A}
\gamma(0)=\dot\gamma(0)=\mathbf 0,
\quad \xi(0)=\pmt{1\\0},\quad
\eta(0)=\frac{250}{27}\mb e_2.
\end{equation}
By \eqref{eq:alpha}, it holds that
$$
\frac{d}{d\tau} \log|[\ddot \gamma,\gamma^{(3)}]|
=\frac{[\ddot \gamma,\gamma^{(4)}]}{[\ddot \gamma,\gamma^{(3)}]}
=\beta_2=\frac{d}{d\tau}\log|18+25 \tau^2 h|,
$$
and there exists a positive constant $C$ such that
$$
[\ddot \gamma,\gamma^{(3)}]=C(18+25 \tau^2 h).
$$
Then \eqref{eq:initial_A} yields that
$C=125/243$, and we get 
\begin{equation}\label{eq:d23A}
[\dot \gamma,\ddot \gamma]=125\tau^2/27,
\end{equation}
which implies that $\tau=0$ is a $3/2$-cusp and $\tau$ is the
$3/5$-arclength parameter of $\gamma$.
Moreover, \eqref{eq:d23A} yields that 
$$ 
[\dot \gamma,\gamma^{(3)}]=
\frac{d}{d\tau} [\dot \gamma,\ddot \gamma]
=\frac{250\tau}{27}
$$ 
and
$$ 
[\dot \gamma,\gamma^{(4)}]=
\beta_1 [\dot \gamma,\ddot \gamma]+
\beta_2 [\dot \gamma,\gamma^{(3)}]
=-\left(\frac{5}{3}\right)^5 \tau^2 h.
$$ 
Using these identities, one can easily check that
$(s_A)^2\kappa_A$ is equal to $(4/25)+\tau^2h$.
\qed

\begin{figure}[htb]
 \begin{center}
        \includegraphics[height=3.0cm]{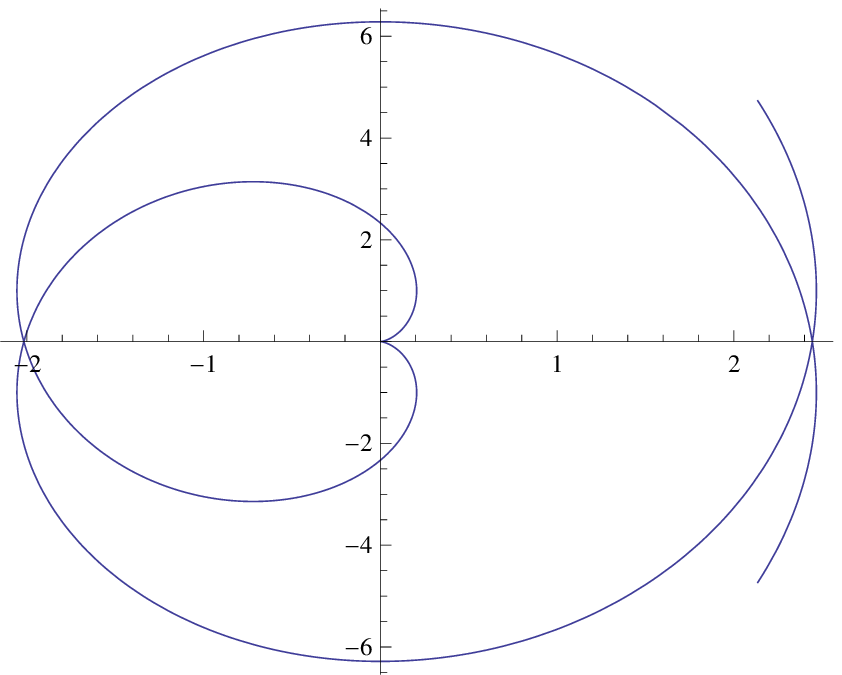}
\qquad \qquad        \includegraphics[height=3.0cm]{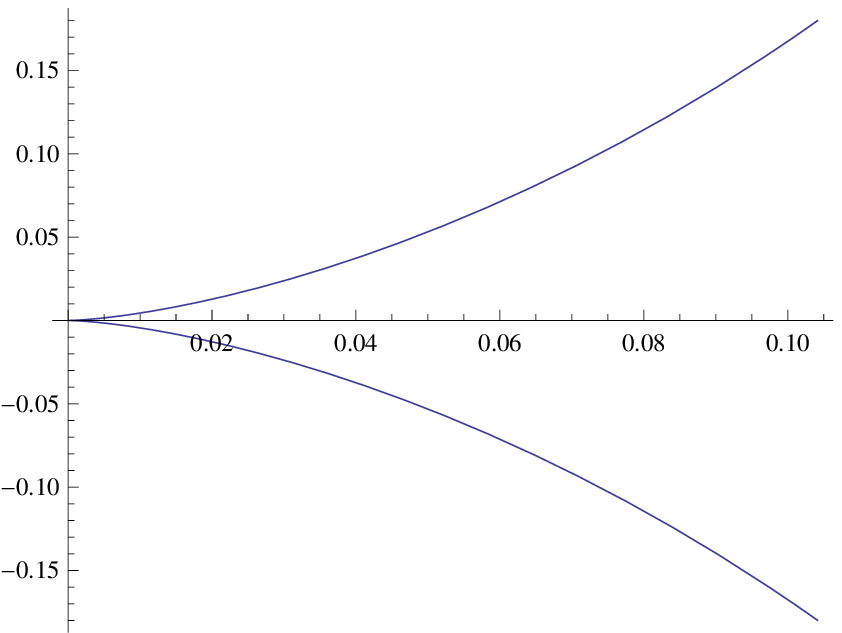}
\caption{The $3/2$-cusps $\sigma_1$ (left) and $\sigma_{-1}$ (right) }
\label{fig:affine}
\end{center}
\end{figure}

\begin{example}
The normalized affine curvature function of
the cusp 
$
\gamma_0(t)={}^t\!\pmt{ t^2 , t^3}
$
is identically equal to $4/25$.
The parameter $t$ is 
proportional to the $3/5$-arclength parameter. 
It is interesting to consider the family of curves
$\sigma_c(\tau)$ with $3/5$-arclength parameter
$\tau$ whose normalized affine curvature function is
equal to $(4/25)+c\tau^2$.
When $c=0$, then $\sigma_0$ is equal to the cusp $\gamma_0$.
By solving the ordinary differential equation \eqref{ode:A1}
with  the initial data \eqref{eq:initial_A}, we get the figure of
$\sigma_c(\tau)$ for $c=\pm 1$ (see Figure \ref{fig:affine}).
\end{example}

\begin{example}
The cycloid \eqref{eq:cyc}
has a positive cusp $t=0$ with 
positive affine cuspidal curvature $\mu_A=36a^{-4/5}$.
(The sign of the affine cuspidal curvature
is independent of the signature of cusps.)
\end{example}

\begin{example}
The hyperbolic cycloid
$$
\gamma(t)=a \pmt{t-\sinh t\\ -1+\cosh t}\qquad (a>0)
$$
has a positive cusp $t=0$ with 
negative affine cuspidal curvature $\mu_A=-36a^{-4/5}$.
\end{example}


To finish this section, we generalize 
Theorem \ref{thm:k_A} for cusps
in an arbitrary  2-manifold with an
equi-affine structure (cf. \cite{IS}):
Let $M^2$ be an oriented $2$-manifold
with an affine connection $D$ (which may not be
torsion free).
If there exists a non-vanishing parallel positively
oriented $2$-form $\Omega$ defined on $M^2$
with respect to $D$,
then the triple $(M^2,D,\Omega)$ is called
an {\it equi-affine 2-manifold}.
(In the usual definition of equi-affine structure,
$D$ is a torsion free connection having symmetric Ricci tensor. However,
we do not assume this here, since it is not
needed for the following discussions.)
For the sake of simplicity, we set
$$
[v,w]:=\Omega(v,w)\qquad (v,w\in T_pM^2,\,\, p\in M^2).
$$
Let $\gamma(t)$ ($|t|<\delta$) be a regular 
curve in $M^2$.
Then the affine curvature function 
is defined by exactly the same formula \eqref{eq:def_k}, 
where
$$
\ddot \gamma(t)=D_t\dot\gamma(t),\quad
\gamma^{(3)}=D_t\ddot\gamma,\quad
\gamma^{(4)}=D_t\gamma^{(3)}.
$$
If $(M^2,g)$ is an orientable Riemannian $2$-manifold,
then the area element $\Omega_g$
is globally defined on $M^2$ 
which is a parallel $2$-form with respect to the Levi-Civta connection. 
Then $(M^2,\nabla,\Omega_g)$ gives
a typical example of an equi-affine 2-manifold.

Let $\gamma(t)$ ($|t|<\delta$) be a $3/2$-cusp 
at $t=0$ on $M^2$. We define the {\it affine cuspidal
curvature} $\mu_A$ at $t=0$ by \eqref{eq:mu_A},
as a generalization of the case of $\R^2$.
Then the following assertion holds, which is a 
generalization of Theorem \ref{thm:k_A}.

\begin{theorem}\label{thm:g_A}
Let $\gamma(t)$ $(|t|<\delta)$ be a $3/2$-cusp 
at $t=0$ in an equi-affine $2$-manifold $(M^2,D,\Omega)$.
Then the same assertion as in Theorem \ref{thm:k_A} holds.
\end{theorem}

By Lemma \ref{lem:parallel},
we have the expression
$$
\dot \gamma(t)
=P_{2}t+\frac{P_{3}}2 t^2
+t^3\Gamma(t),
$$
where $P_j:=P_j(t)$ ($j=1,2,3,...$)
is a parallel vector field along $\gamma(t)$
such that $P_j(0)=\gamma^{(j)}(0)\in T_{\gamma(0)}M^2$, and
$\Gamma(t)$ is a $C^\infty$-vector field along $\gamma$.
Using this expression, 
like as in the proof of Theorem \ref{thm:g_g},
we can prove Corollary \ref{cor:k_A}.
The remainder of the proof is exactly the same as that of 
Theorem \ref{thm:k_A}.
Like as in the Euclidean case
 (cf. \eqref{eq:map_A}), the map 
\begin{equation}\label{eq:map_A}
\mathcal F_A:S^*_0(\R,M^2)
\ni \gamma(t)\mapsto s_A(t)^2\kappa_A(t)
\in C^\infty_0(\R)
\end{equation}
is defined, namely $\mathcal F_A(\gamma)$ is the
normalized affine curvature function of the $3/2$-cusp $\gamma$
with respect to the equi-affine structure $(D,\Omega)$.
The following assertion follows immediately.

\begin{corollary}\label{cor:image_A}
The image of the map $\mathcal F_A$ coincides with 
the subset
$$
\Sigma_A:=
\left\{f\in C^\infty_0(\R)\,;\, f(0)=\frac{4}{25},\,\, \dot f(0)=0
\right\}.
$$
\end{corollary}

\section{Affine geometry of generic inflection points}

A point $t=c$ of a regular curve $\gamma(t)$ in the
affine plane $\R^2$ is called an {\it inflection point} if
$[\dot \gamma(c),\ddot \gamma(c)]$ vanishes. 
It is well-known that there is a 
duality between cusps and inflection points
(cf. \cite{SUY}).
An inflection point $t=c$
is called {\it generic} if 
$[\dot \gamma(c),\gamma^{(3)}(c)]$ does not vanish.
Fabricius-Bjerre \cite{FB} pointed out that
$\kappa_A$ diverges to $-\infty$ for generic
inflection points.
More generally, $\lim_{c\to 0}\inf_{|t|<c}\kappa_A(t)$
diverges to $-\infty$ 
for arbitrary inflection points (see \cite{U}). 

\begin{definition}
A generic inflection point $t=c$ is called {\it positive}
(resp. {\it negative}) if $[\dot \gamma(c),\gamma^{(3)}(c)]$
is positive (resp. negative).
\end{definition}

This signature of inflection points is invariant under
the orientation preserving diffeomorphisms of $\R^2$
and is also invariant under
the choice of an orientation of curves.
In this section, we define a new 
affine invariant called
the {\it affine inflectional curvature} 
at an inflection point $t=c$ by
\begin{equation}\label{eq:mu_I}
\mu_I:= \epsilon_I \frac{[\dot \gamma(c),\gamma^{(4)}(c)]
-6 [\ddot \gamma(c),\gamma^{(3)}(c)]}{
[\dot \gamma(c),\gamma^{(3)}(c)]^{5/4}},
\end{equation}
where $\epsilon_I=\op{sgn}[\dot \gamma(c),\dddot \gamma(c)]$ 
is the signature of the inflection point.
This is invariant under equi-affine transformations, but
changes sign if we reverse the orientation 
of the curve $\gamma$. 
As an analogue of Theorem \ref{thm:main},
we get the following assertion:

\begin{theorem}\label{thm:main_I}
Let $\gamma_1(t)$ and $\gamma_2(t)$ be two positive
generic inflection points at $t=0$ in the
affine plane $\R^2$. 
Then $\gamma_1$ and $\gamma_2$ have 
the same inflectional curvature
if and only if there exists an 
orientation preserving equi-affine transformation $T$
of $\R^2$ and a parametrization $u=u(t)$ near $t=0$
such that $du/dt>0$ and
$|T\circ \gamma_2\circ u(t)-\gamma_1(t)|$
has order higher than $t^4$.
\end{theorem}

This assertion immediately follows from
the following proposition.

\begin{proposition}\label{prop:equi2}
Let $\gamma(t)$ be a regular curve in $\R^2$
having 
a positive generic inflection point at $t=0$.
Then there exist an orientation preserving equi-affine transformation $T$
and
a coordinate change $t=t(u)$
such that
\begin{equation}
T\circ \gamma\circ t(u)=\pmt{u\\ u^3 +\sqrt[4]{6} \mu_Iu^4/4} 
+o(u^4),
\end{equation}
where 
$o(u^4)$ is a term of order higherthan $u^4$.
\end{proposition}

Since the proof of the proposition is similar to that of
Proposition \ref{prop:equi2}, we omit it here.
We also get the following assertion:

\begin{theorem}\label{thm:k_I}
Let $t=0$ be a generic inflection point
of a regular curve $\gamma(t)$ in $\R^2$.
Then  $\tau=\op{sgn}(t)(s_A)^{3/4}$
can be taken to be a coordinate of $\gamma$ at $t=0$
{$($called the \em  $3/4$-arclength parameter)}
and $f:=(s_A)^2\kappa_A$ is a $C^\infty$-function 
of $t$ $($and $\tau)$. Moreover,
it satisfies 
\eqref{eq:hidden2},
and
\begin{align}\label{eq:lim_I}
&\lim_{t\to 0}(s_A)^2\kappa_A(=f(0))=-\frac{5}{16},\\
\label{eq:lim_I2}
&\lim_{t\to 0} \frac{(s_A)^2\kappa_A+{5}/{16}}{\tau}
=-\frac{3 \sqrt[4]{3}}{28 \sqrt{2}}\mu_I.
\end{align}
Furthermore, if $t=\tau$, then \eqref{eq:hidden2} is reduced
to the relation
\begin{equation}\label{eq:hidden}
32\left(\frac{df(0)}{d\tau}\right)^2
+9\frac{d^2f(0)}{d\tau^2}=0.
\end{equation}
Conversely, if we take a $C^\infty$-function
$f(\tau)$ satisfying $f(0)=-5/16$
and \eqref{eq:hidden},
then there exists a generic inflection point
whose normalized affine curvature function is
$f(\tau)$ with respect to the $3/4$-arclength
parameter.
\end{theorem}

\begin{proof}
Let $t=0$ be a generic inflection point
of a plane curve $\gamma(t)$.
Since orientation reversing eqi-affine transformations 
preserve $\kappa_A$ and $\mu_I$,
we may assume that $t=0$ is a positive inflection point.
As an analogue of Proposition \ref{prop:k_A},
$\tau=\op{sgn}(t)(s_A)^{3/4}$ can be
taken to be a coordinate of $\gamma$ at $t=0$.
(As in Theorem \ref{thm:k_I},
$\tau$ is called the $3/4$-arclength
parameter.)
Also, we can prove that $f(t):=(s_A)^2\kappa_A$ is
a $C^\infty$-function of $t$.

The following identities hold (cf. \eqref{eq:exponent})
\begin{equation}\label{eq:i12}
\dot \gamma=\frac43\tau^{1/3}\gamma', \qquad
\ddot \gamma=
\frac{4 \left(\gamma'+4 \tau^{4/3}\gamma'' \right)}{9 \tau^{2/3}},
\end{equation}
where we use the notations as in 
\eqref{eq:dot} and \eqref{eq:prime},
namely, 
the dot (resp. the prime)
means the derivative with respect to 
the $3/4$-arclength parameter
$\tau$ (resp. the affine arclength parameter $s_A$).
Using \eqref{eq:i12}, one can easily check that
$t$ is a $3/4$-arclength parameter 
(of a positive inflection point)
if and only if
$[\dot \gamma,\ddot \gamma]=64t/27$.
By \eqref{eq:i12}, it holds that
$$
\gamma^{(3)}=
\frac{8 \left(-\gamma'+6\tau^{4/3}\gamma'' 
+8\tau^{8/3}\gamma^{[3]} 
\right)}{27 \tau^{5/3}}.
$$
Substituting $\gamma^{[3]}=-\tau^{-8/3}f\gamma'$,
we get
$$
\gamma^{(3)}=
-\frac{8 \left((8 f+1) \gamma'-6 \tau^{4/3} 
\gamma''\right)}{27 \tau^{5/3}}.
$$
Differentiating it by using $\gamma^{[3]}=-\tau^{-8/3}f\gamma'$
again, we have
$$
\gamma^{(4)}=\frac{8}{81\tau^{8/3}}
\left((-24 \tau \dot f+16 f+5) 
\gamma '-2 \tau^{4/3} (16 f+5) \gamma ''\right).
$$
Since $[\gamma',\gamma'']=1$, it holds that 
\begin{equation}\label{eq:3}
[\dot\gamma,\ddot \gamma]=\frac{64}{27}\tau,\quad
[\dot\gamma,\gamma^{(3)}]=\frac{64}{27},\quad
[\ddot\gamma,\gamma^{(3)}]
=\frac{64 (5 + 16 c + 16 \tau g(\tau))}{243 \tau},
\end{equation}
where we set
$
(s_A)^2\kappa_A=c+\tau g(\tau) \,\, (g\in C^\infty_0(\R)).
$
Since $[\ddot\gamma,\gamma^{(3)}]\in C^\infty_0(\R)$,
we have that $c=-5/16$. 
Moreover, since $[\dot \gamma,\gamma^{(3)}]$
is a constant function, 
$[\dot \gamma,\gamma^{(4)}]$ and $-[\ddot \gamma,\gamma^{(3)}]$
are both equal to
$-{2^{10}g(\tau)}/{3^5}$.
Using these relations
and \eqref{eq:mu_I}, 
one can easily check that
$$
-\frac{3 \sqrt[4]{3}}{28\sqrt{2}}\mu_A=g(0)
\left(=\lim_{\tau\to 0} \frac{(s_A)^2 \kappa_A+5/16}{\tau}\right),
$$
which proves the identity \eqref{eq:lim_I2}. 
On the other hand, it holds that
$$
[\gamma^{(3)},\gamma^{(4)}]=
\frac{2^{10} \left(9 \dot g(\tau)+16 g(\tau)^2\right)}{3^7 \tau}.
$$
Since $[\gamma^{(3)},\gamma^{(4)}]\in C^\infty_0(\R)$,
we can conclude 
$9 \dot g(0)+16 g(0)^2=0$,
which is equivalent to the condition 
\eqref{eq:hidden}.
In this situation,
there exists a smooth function $h(\tau)$
near $\tau=0$ such that
\begin{equation}\label{eq:g0}
9 \dot g(\tau)+16 g(\tau)^2=\tau h(\tau).
\end{equation}
Since $\dot \gamma$ and $\gamma^{(3)}$ are linearly independent,
we can write
\begin{equation}\label{eq:ode_I20}
\ddot \gamma=a_{11}\dot \gamma +a_{12}\gamma^{(3)},\quad
\gamma^{(4)}=a_{21}\dot \gamma+a_{22}\gamma^{(3)}.
\end{equation}
By \eqref{eq:g0}, we have
\begin{align}\label{eq:a}
a_{11}&=\frac{[\ddot \gamma,\gamma^{(3)}]}{[\dot \gamma,\gamma^{(3)}]}
=\frac{16 g(\tau)}9,\quad a_{12}
=\frac{[\dot \gamma,\ddot\gamma]}{[\dot \gamma,\gamma^{(3)}]}
=\tau, \\ \label{eq:a2}
\quad a_{21}&=
-
\frac{[\gamma^{(3)},\gamma^{(4)}]}{[\dot \gamma,\gamma^{(3)}]}=
-\frac{16 h(\tau)}{81}, 
\quad a_{22}=
\frac{[\dot\gamma,\gamma^{(4)}]}{[\dot \gamma,\gamma^{(3)}]}
=
-\frac{16 g(\tau)}{9}.
\end{align}
We now take a function 
$f\in C^\infty_0(\R)$
satisfying \eqref{eq:hidden}.
Let $g(\tau)$ 
be the function defined by
$g(\tau)=\tau f(\tau)$,
and $h(\tau)$ be 
the function given by
\eqref{eq:g0}.
We consider the 
ordinary differential equation
\begin{equation}\label{eq:ode_I2}
\dot \gamma=\xi,\quad \dot \xi=a_{11}\xi+a_{12}\eta,\quad
\dot\eta=a_{21}\xi+a_{22}\eta,
\end{equation}
where $a_{ij}$ ($i,j=1,2$)
are defined by \eqref{eq:a} and \eqref{eq:a2}. 
Then there exists a solution $\gamma(\tau)$ 
of the ordinary differential equation \eqref{eq:ode_I2}
satisfying the initial conditions
$$
\gamma(0)=\mathbf 0,
\quad \xi(0)=\pmt{1 \\ 0},\quad
\eta(0)=\frac{64}{27}\mb e_2.
$$
By \eqref{eq:ode_I2} and \eqref{eq:a},
it can be easily checked that 
$
\gamma^{(3)}=\ddot \xi=\eta.
$
Then $\phi:=[\dot \gamma,\ddot \gamma]$
satisfies
$$
\phi=[\xi,\dot \xi]
=a_{12}[\xi, \eta]=
a_{12} [\dot \gamma,\gamma^{(3)}]=a_{12} \dot \phi,
$$
which  yields that
$$
\frac{d}{d\tau}\log |\phi|
=\frac{1}{a_{12}}=\frac{1}{\tau}
=\frac{d}{d\tau}\log |\tau|.
$$
Thus there exists a positive constant $C$
such that $\phi=C\tau$.
Since 
$$
\dot\phi(0)
=[\dot \gamma(0),\gamma^{(3)}(0)]
=[\xi(0),\eta(0)]=\frac{64}{27},
$$
we have $C=64/27$.
Then 
\begin{equation}\label{eq:d130i}
[\dot \gamma(\tau),\gamma^{(3)}(\tau)]=\frac{64}{27},
\end{equation}
which implies that $\tau=0$ is a generic inflection point
and $\tau$ is the $3/4$-arclength parameter.
Moreover, we can get the relation
$[\dot \gamma,\gamma^{(4)}]=
-{2^{10}g(\tau)}/{3^5}$
from the identity 
$\gamma^{(4)}=a_{21}\dot \gamma+a_{22}\gamma^{(3)}$.
Since $[\dot \gamma,\gamma^{(3)}]$ is a constant,
\begin{equation}\label{eq:d23d14}
[\ddot \gamma,\gamma^{(3)}]=
-[\dot \gamma,\gamma^{(4)}]
\end{equation}
holds. Using these relations, it can be easily checked that
$(s_A)^2\kappa_A$ is equal to $-({5}/{16})+\tau g(\tau)$, which proves
 last assertion of Theorem \ref{thm:k_I}.

Finally, using the initial parameter $t$ of $\gamma(t)$,
we have that
$$
32 (\dot f)^2+9 \ddot f=
32 \left(\frac{dt}{d\tau}\right)^2
f_t^2+9\frac{d^2t}{d\tau^2}f_t
+9 \left(\frac{dt}{d\tau}\right)^2f_{tt},
$$
where
$f_t:=df/dt$
and $f_{tt}:=d^2f/dt^2$.
Since
$$
\frac{d^2 t/d\tau^2}{(dt/d\tau)^2}
=-\frac{d}{d\tau}\left(\frac{1}{dt/d\tau}\right)
=-\frac{d t}{d\tau}\frac{d}{dt}\left(\frac{d\tau}{dt}\right)
=-\frac{d^2\tau/dt^2}{d\tau/dt},
$$
\eqref{eq:hidden2} is reduced to the relation
$$
32 f_t(0)^2-9\frac{d^2\tau(0)/dt^2}{d\tau(0)/dt}f_t(0)
+9 f_{tt}(0)=0.
$$
Since $\tau=(s_A)^{3/4}$, one can prove the 
following identity using L'Hospital's rule
\begin{align*}
\frac{d\tau(0)}{dt}&=
\left(\frac{3}4\right)^{3/4} 
[\gamma_t(0),\gamma_{ttt}(0)]^{1/4}, \\
\frac{d^2\tau(0)}{dt^2}&=
\left(\frac{3}4\right)^{3/4}\frac{
 [\gamma_t(0),\gamma_{tttt}(0)]+[\gamma_{tt}(0),
\gamma_{ttt}(0)]}
{7[\gamma_t(0),\gamma_{ttt}(0)]^{3/4}},
\end{align*}
which yield the identity \eqref{eq:hidden2},
where $\gamma_t:=d\gamma/dt$
and $\gamma_{tt}:=d^2\gamma/dt^2$ etc.
\end{proof}
\begin{example}
The normalized affine curvature function of
the inflection point 
$
\gamma(t)={}^t\!\pmt{ at , at^3} \,\,(a>0)
$
is identically equal to $-5/16$. 
The parameter $t$ is 
proportional to the $3/4$-arclength parameter. 
\end{example}

\begin{example}
The {\it skew-cycloid}
$$
\gamma(t)=a\pmt{t-\sin t, -t+\cos t}\qquad (a>0)
$$
has a positive inflection point
at $t=0$ 
of inflectional curvature $\mu_I=-{6}/{\sqrt{a}}$.
If one reverse the orientation, $\gamma(-t)$
also has a positive inflection point
with the inflectional curvature $\mu_I={6}/{\sqrt{a}}$.
\end{example}


Like as in the case of cusps,
the following assertion can be proved
by modifying the proof of 
Theorem \ref{thm:k_I}:

\begin{theorem}\label{thm:g_I}
Let $\gamma(t)$ $(|t|<\delta)$ be a generic
inflection point at $t=0$ in an equi-affine 
$2$-manifold $(M^2,D,\Omega)$.
Then the same assertion as in Theorem \ref{thm:k_I}
holds.
\end{theorem}

We denote by $I^*_0(\R,M^2)$ the 
set of germs of $C^\infty$-maps $\gamma(t)$ 
defined on an open interval containing $t=0$
into an equi-affine $2$-manifold $(M^2,D,\Omega)$
which gives a generic inflection point at $t=0$.
Then the map 
\begin{equation}\label{eq:map_I}
\mathcal F_I:I^*_0(\R,M^2)
\ni \gamma(t)\mapsto s_A(t)^2\kappa_A(t)
\in C^\infty_0(\R)
\end{equation}
is defined, namely $\mathcal F_I(\gamma)$ is the
normalized affine curvature function of the 
generic inflection point $\gamma$.

\begin{corollary}\label{cor:image_I}
The image of the map $\mathcal F_I$ coincides with 
the subset
{
\begin{align*}
\Sigma_I
&:=\left\{
f\in C^\infty_0(\R)\,;\, f(0)=-\frac{5}{16},\,\,\dot f(0)\ne 0
\right\}
\\
&\phantom{******}\cup 
\left\{f\in C^\infty_0(\R)\,;\, 
f(0)=-\frac{5}{16},\,\,\dot f(0)=\ddot f(0)=0
\right\}.
\end{align*}}
\end{corollary}

\begin{proof}
Obviously the image of $\mathcal F_I$ is contained
in $\Sigma_I$.
We fix a function $f\in C^\infty_0(\R)$
such that $f(0)=-5/16$. 
If $\dot f(0)=\ddot f(0)=0$,
then $f$ satisfies \eqref{eq:hidden}.
On the other hand, if $\dot f(0)\ne 0$,
then there exists a new parametrization 
$\tau:=t+c t^2$ such that $f(\tau):=f\circ t(\tau)$
satisfies \eqref{eq:hidden},  
by adjusting the constant $c$.
Thus, for these two cases, 
Theorem \ref{thm:g_I} yields that
there exist $\delta>0$ and
a regular curve 
$\gamma:(-\delta,\delta)\to (M^2,g)$
having  an inflection point $\tau=0$
such that $s_A(\tau)^2\kappa_A(\tau)=f(\tau)$
and $\tau$ is the $3/4$-arclength parameter.
Since the relation
$s_A(\tau)^2\kappa_A(\tau)=f(\tau)$
is independent of the choice of parameters,
we get the assertion.
\end{proof}

\appendix

\section*{Appendix. A division lemma}

\begin{lemma}\label{lem:fractional}
Let $\phi(t)$ be a $C^\infty$-function at $t=0$,
and $\alpha$ a positive real number.
Then the function defined by
$$
f(t):=\frac{\Phi(t)}{\op{sgn}(t)|t|^{1+\alpha}}
\qquad (\Phi(t):=\int_0^t |u|^\alpha \phi(u) du)
$$
is a $C^\infty$-function at $t=0$, namely
$f\in C^\infty_0(\R)$.
Moreover, it holds that
$$
f(0)=\frac{\phi(0)}{1+\alpha}.
\leqno (A.1)
$$
\end{lemma}

\begin{proof}
In fact, we have that
$$
\Phi(t)=\int_0^1 \frac{d\Phi(tu)}{du} du
=\int_0^1 t \dot \Phi(tu) du
=\op{sgn}(t)|t|^{1+\alpha} 
\left(\int_0^1 |u|^\alpha \phi(tu) du\right).
$$
Since $\alpha>0$, it follows that
$\int_0^1 |u|^\alpha \phi(tu) du$ 
is a $C^\infty$-function at $t=0$.
By L'Hospital's rule, we have (A.1).
\end{proof}

If $\alpha=0$, then $f(t)=\Phi(t)/t$,
and the lemma reduces to the classical division lemma.


\begin{thebibliography}{KR}
\bibitem{FB}
Fr. Fabricius-Bjerre,
{\it On a conjecture of Bol},
Math. Scand. 40 (1977), 194--196. 

\bibitem{GS}
P. Giblin and G. Sapiro,
{\it Affine-Invariant Distances, Envelopes and
Symmetry Sets},
Geometriae Dedicata 71 (1998) 237--261.

\bibitem{IS}
S. Izumiya and T. Sano,
{\it Generic affine differential geometry of plane curves}, 
Proc. Math. Soc. Edinburgh. 41 (1998) 315--324. 

\bibitem{NS}
K. Nomizu and T. Sasaki, 
{\it Affine Differential Geometry}, 
Cambridge University Press, Cambridge, 1994.

\bibitem{SUY}
K. Saji, M. Umehara and K. Yamada, 
{\it 
The duality between singular points and
inflection points on wave fronts},
Osaka J. Math. 47 (2010), 591--607.

\bibitem{U0}
M. Umehara, {\it Differential geometry on surfaces with singularities}
 in The World of Singularities
(ed. H. Arai, T. Sunada and K. Ueno) Nippon-Hyoron-sha Co., Ltd. 
(2005), 50--64, (Japanese).

\bibitem{U}
M. Umehara,
{\it A simplification of the proof
of Bol's conjecture on sextactic points},
Proc. Japan Acad. Ser. A  87 (2011), 10-12. 

\end{thebibliography}
\end{document}